\documentclass[letterpaper, 11 pt]{article}  % Comment this line out if you need a4paper
\usepackage[margin=1in]{geometry}
\usepackage{graphicx}	
\usepackage{pdfpages}		
\usepackage{amssymb}
\usepackage{empheq}
\usepackage{color}
\usepackage{mathabx}
\usepackage{adjustbox}
\usepackage{centernot}
\usepackage{tabularx}
\usepackage{arydshln}
\usepackage{mathtools}
\usepackage{cases}
\usepackage{mathrsfs}  
\usepackage{csquotes}
\usepackage{epstopdf}
\usepackage{mathtools}
\usepackage{amsmath}
\usepackage{amsfonts}
\usepackage{nicefrac}
\usepackage{hyperref}
\usepackage{nccmath}
\usepackage{pst-node}
\usepackage{csquotes}
\usepackage{arydshln}
\usepackage{multirow}
\usepackage{blkarray}
\usepackage[ruled]{algorithm2e}
\usepackage{algcompatible}
\usepackage{soul}
\allowdisplaybreaks

\newcommand{\bb}{\mathbf}
\newcommand{\mc}{\mathcal}
\newcommand{\nn}{\nonumber}
\usepackage{amssymb}
\newtheorem{theorem}{\textbf{Theorem}}

\newtheorem{lemma}{\textbf{Lemma}}
\newtheorem{prop}{\textbf{Proposition}}

\newcommand{\tss}{\textsuperscript}
\newcommand{\tsb}{\textsubscript}
\newcommand{\te}{\text}

\DeclareMathOperator*{\minimize}{minimize}

\title{\LARGE \bf
Co-Design of Delays and Sparse Controllers for Bandwidth-Constrained Cyber-Physical Systems}
\author{Nandini Negi\tss{1,2} and Aranya Chakrabortty\tss{1,3}\\
\tss{1}Electrical $\&$ Computer Engineering, North Carolina State University\\
Email : \tss{2}{nnegi@ncsu.edu}, \tss{3}{aranya.chakrabortty@ncsu.edu}
\thanks{The research presented in this paper was partly supported by the US National Science Foundation under grant ECCS 1509137.}
}
\date{\vspace{-1cm}}
\begin{document}
\setstcolor{cyan}
\maketitle
\thispagestyle{empty}
\pagestyle{empty}
\begin{abstract}
We address the problem of sparsity-promoting optimal control of cyber-physical systems with feedback delays. The delays are categorized into two classes - namely, intra-layer delay, and inter-layer delay between the cyber and the physical layers. Our objective is to minimize the $\mc{H}_2$-norm of the closed-loop system by designing an optimal combination of these two delays along with a sparse state-feedback controller, while respecting a given bandwidth constraint. We propose a two-loop optimization algorithm for this. The inner loop, based on alternating directions method of multipliers (ADMM), handles the conflicting directions of decreasing $\mc{H}_2$-norm and increasing sparsity of the controller. The outer loop comprises of semidefinite program (SDP)-based relaxations of non-convex inequalities necessary for stable co-design of the delays with the controller. We illustrate this algorithm using simulations that highlight various aspects of how delays and sparsity impact the stability and $\mc{H}_2$-performance of a LTI system.

\end{abstract}

\section{Introduction}

\par \noindent In recent years, sparsity-promoting optimal control has emerged as a key tool for enabling economical control of large-scale cyber-physical systems (CPSs). 
% These designs aim to minimize the number of communication links needed for control of the CPS while guaranteeing a desired closed-loop performance. Several variants of this design have been proposed using optimization algorithms 
such as ADMM \cite{mihailo}, LASSO \cite{zico}, GraSP \cite{jsac}, and PALM \cite{linpalm}. An extension of these results to LTI systems with communication delays has been reported in \cite{negi2018}. Since most real-world CPSs operate under stringent constraints for bandwidth, stability and closed-loop performance in the presence of delays are important requirements for these controllers \cite{delay-hespanha}. Accordingly, the algorithm in \cite{negi2018} derives convex relaxations of bilinear matrix inequalities to design a sparse controller, while guaranteeing closed-loop stability under a constant delay.
\par In this paper, we extend the design in \cite{negi2018} one step further by considering the delays themselves as {\it design variables}. Our formulation is motivated by modern CPS communication technologies such as software-defined networking (SDN) and cloud computing that offer flexibility to network operators in choosing delays in communication links. We consider two kinds of delays - namely (1) \textit{inter-layer} delay that arises in the local-area network (LAN) connecting the sensors in the physical layer to the computational units in the cyber layer, and (2) intra-layer delay that arises in the SDN connecting the computational units spread across the cyber-layer. Our goal is to co-design these two delays with a sparse feedback controller so that the $\mc{H}_2$-norm of the closed-loop system is minimized, while ensuring that both delays are greater than or equal to their individual lower bounds that arise from the cost of the network bandwidth. The main contribution of this paper is to develop a hierarchical optimization algorithm that provides a guided solution for this co-design. The outer loop designs the two delays and finds a corresponding stabilizing controller by sequentially relaxing the non-linear matrix equations required for the co-design. The inner loop sparsifies this controller while minimizing the closed-loop $\mc{H}_2$-norm. Our results show that depending on the plant dynamics, the relative magnitudes of the two delays for achieving the optimal $\mc{H}_2$-norm can be notably different.
% \par Compared to \cite{negi2018}, this extension, however, is not straightforward. This is because \cite{negi2018} considered both the intra-layer and inter-layer delays to be equal and fixed. Considering them to be design variables comes at the price of adding an extra layer of complexity as now one must co-optimize both the controller and the delays, which are all coupled to each other through complex implicit relationships arising from stability \cite{stability}, $\mc{H}_2$-performance, and bandwidth constraints. To handle these dependencies, we design our algorithm with two hierarchical loops.  
\par Note that our problem is fundamentally different from the conventional bandwidth allocation and delay assignment problems commonly addressed in the networking literature \cite{net1}, \cite{net2}, where the utility functions to be optimized are static objectives. Our goal, in contrast, is to design a bandwidth allocation mechanism that minimizes the $\mc{H}_2$-norm of a CPS over a sparse state-feedback controller. 
% The objective function for our design, therefore, addresses the inherent {\it dynamics} of the plant and sensitivity of its stability as a function of the two delays, as well as the sparse structure of the control.
We illustrate the effectiveness of our algorithm using simulations that highlight the impacts of delays and sparsity on $\mc{H}_2$-performance. 
\par The rest of the paper is organized as follows. Section \ref{sec:problemformulation} states the problem formulation followed by Section \ref{sec:proposedcodesign} that describes the proposed co-design of the delays. Section \ref{sec:problemsetupinadmm} introduces the two-loop algorithm to solve the problem followed by simulations in Section \ref{sec:simulations}, and conclusion in Section VI. The proofs of all lemmas, theorems and propositions are listed in the Appendix unless stated otherwise.

\section{Problem Formulation}
\label{sec:problemformulation}
\subsection{State Feedback with Communication Delays}
\begin{figure}[hbtp!]
\centering
\includegraphics[scale=0.5]{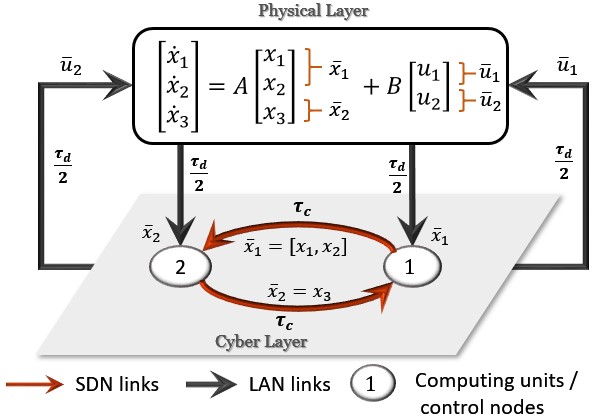}
\caption{Sample CPS schematic showing physical and cyber layers with the associated delays.}
\label{fig:system}
\vspace{-5mm}
\end{figure}

\noindent Consider a LTI system with the following dynamics:
\begin{equation}
\dot{x}(t)=Ax (t)+ B u(t) + B_w w(t), \label{delay-free}
\end{equation}
where $x \in \mathbb{R}^{n}$ is the state, $u \in \mathbb{R}^{m}$ is the control, and $w \in \mathbb{R}^r$ is the exogenous input, with the corresponding matrices $A \in \mathbb{R}^{n\times n}$, $B\in\mathbb{R}^{n\times m}$ and $B_w\in \mathbb{R}^{n\times r}$. We design a state-feedback controller, ideally represented as $u(t)=-Kx(t)$. However, due to limited bandwidth availability, the controller includes finite delays in the feedback. The CPS model that we consider is described as follows.
\begin{enumerate}
    \item There are $p$ sensors and actuators in the physical layer and the state vector $x(t)$ is correspondingly divided into $p$ non-overlapping parts $\bar{x}_1(t),\ldots,\bar{x}_p(t)$, where $\bar{x}_i$ is measured by the $i$-th sensor.
    \item There are $p$ computing units or control nodes located in a virtual cloud network. The $i$-th sub-state $\bar{x}_i(t)$ is transmitted to the $i$-th control node through LAN with incident delay $\nicefrac{\tau_d}{2}$.
    \item Inside the cloud, also referred to as the \textit{cyber layer}, the control nodes share their individual sub-states $\bar{x}_i(t)$ with each other over an SDN with delay $\tau_c$. Each control node $i$ calculates a portion of the control input vector denoted as $\bar{u}_i(t)\in\mathbb{R}^{m_i}$, where $ \sum_{i=1}^p m_i= m$.
    \item The calculated control inputs are transmitted back to the physical layer with $\nicefrac{\tau_d}{2}$ delay. 
\end{enumerate}
A schematic of this CPS with $n=3$, $m=2$ and $p=2$ is shown in Fig. \ref{fig:system}. Denoting $\tau_o=\tau_d+\tau_c$, the control input can be expressed as:
\begin{align}
&u(t)  = - \underbrace{(K\circ \mc{I}_d)}_{K_d} x(t-\tau_d) - \underbrace{(K\circ \mc{I}_o)}_{K_o} x(t-\tau_o), \label{Ueq}
\end{align}
where $\circ$ represents Hadamard product. $\mc{I}_d,\mc{I}_o\in\mathbb{R}^{m\times n}$ are binary matrices such that
\begin{align}
    \mc{I}_d(i,j) =&\begin{cases}
   1, \ \text{If} \ \exists \ q\in\{1,\ldots,p\} :  u_i \in \bar{u}_q \ \text{and} \ x_j \in \bar{x}_q, \ \\
    0, \ \text{otherwise}.
    \end{cases}
\end{align}
and $\mc{I}_o$ is the complement of $\mc{I}_d$. For the system shown in Fig. \ref{fig:system}, $\mc{I}_d$ and $\mc{I}_o$ are:
\begin{equation}
\mc{I}_d=\begin{blockarray}{cccc}
 \BAmulticolumn{1}{c}{} & \BAmulticolumn{2}{c}{\overbrace{}^{\bar{x}_1}} & \BAmulticolumn{1}{c}{\overbrace{}^{\bar{x}_2}}\\
 \begin{block}{c[ccc]}
\bar{u}_1  & 1 & 1 & 0\\
\bar{u}_2  & 0 & 0 & 1\\
    \end{block}
\end{blockarray}, \ \mc{I}_o=\begin{blockarray}{cccc}
 \BAmulticolumn{1}{c}{} & \BAmulticolumn{2}{c}{\overbrace{}^{\bar{x}_1}} & \BAmulticolumn{1}{c}{\overbrace{}^{\bar{x}_2}}\\
 \begin{block}{c[ccc]}
\bar{u}_1  & 0 & 0 & 1\\
\bar{u}_2  & 1 & 1 & 0\\
    \end{block}
\end{blockarray}.
\end{equation}
The closed-loop system of \eqref{delay-free}-\eqref{Ueq} can be written as:
\begin{align}
\dot{x}(t) =& A x(t) - B K_d x(t-\tau_d) - B K_o x(t-\tau_o) + B_w w(t), \nn \\
z(t) =& C x(t) + D u(t) = \begin{bmatrix}
Q^{\nicefrac{1}{2}} \\ 0
\end{bmatrix}, \ D=\begin{bmatrix}
0\\R^{\nicefrac{1}{2}}
\end{bmatrix},  \label{delayed}
\end{align}
where $z(t)$ is the measurable output, $Q \succeq 0$ and $ R\succ 0$. We make the standard assumption that $(A,B)$ and $(A,Q^{\nicefrac{1}{2}})$ are stabilizable and detectable, respectively \cite[Sec. II]{mihailo}. 

\subsection{Problem Setup}
\par \noindent Our goal is to design a $K$ that minimizes the $\mc{H}_2$-norm of the transfer function from $w(t)$ to $z(t)$ for the time-delayed LTI system (\ref{delayed}). In general, 
the $\mc{H}_2$-performance of \eqref{delayed} will be worse than that of the delay-free system \cite[Section 5.6.1]{stability}. Therefore, reducing both the delays $\tau_d$ and $\tau_c$ will improve the $\mc{H}_2$-performance. The trivial solution, of course, would be to use $\tau_d=\tau_o=0$, which is not possible in reality as that would require infinite bandwidth.
\par Let the combined bandwidth of links connecting the physical sensors to the cloud be $W_{cp}$, and that of SDN links inside the cloud be $W_{cc}$. Then, the total cost for renting bandwidth can be written as:
\begin{align}
 S =m_{cp} W_{cp} + m_{cc} W_{cc},
\end{align} 
where $m_{cp}$ and $m_{cc}$ are the respective dollar costs for renting LAN and SDN links. $W_{cp}$ and $W_{cc}$ are divided into the total number of links as described below.
\begin{itemize}
    \item The uplink for carrying $\bar{u}_i(t)$ back to the physical actuator is not needed if the $i$-th block row of $K$ is entirely $0$. Similarly, if the $i$-th block column of $K$ is $0$, then $\bar{x}_i$ is no longer required for calculating any control input, and the corresponding downlink becomes redundant. The uplinks and downlinks together constitute the LAN links. Thus, $W_{cp}$ is effectively divided into the number of non-zero block rows and columns of $K$ denoted by $N_{row}(K)$ and $N_{col}(K)$, respectively. 
    \item $W_{cc}$ is divided into the number of non-zero off-diagonal blocks of $K$ denoted by $N_{off}(K)$.
\end{itemize}
Accordingly, we can write the bandwidth constraint as:
\begin{align}
&S= 2 m_{cp}\left( \frac{N_{row}(K) + N_{col}(K)}{\tau_d}\right) + m_{cc} \left(\frac{N_{off}(K)}{\tau_o-\tau_d}\right) \leq S_b,  \label{bw2}
\end{align} 
where $S_b>0$ is a mandatory budget that is imposed to prevent infinite bandwidth. To minimize the cost of renting the links and bandwidth, we wish to reduce the number of both LAN and SDN links by promoting sparsity in $K$. Our design objectives, therefore, are listed as:
\par \textbf{P1} : Design $\tau_d$, $\tau_o$ and $K$ such that 
\begin{itemize}
\item[$\bullet$] \textit{$\mc{H}_2$-norm of the closed-loop transfer function of \eqref{delayed} from $w(t)$ to $z(t)$, denoted as $J$, is minimized.}
\item[$\bullet$] \textit{The bandwidth cost $S$ satisfies \eqref{bw2} for some given budget $S_b$, which is assumed to be large enough for the problem to be feasible.}
\item[$\bullet$] \textit{Sparsity of $K$ is promoted.}
\end{itemize}
\par \noindent Let $\mathbb{K}$ be the set of all $K$ that stabilize \eqref{delayed} for given delays $\tau_o$ and $\tau_d$. Given $S_b$, \textbf{P1} can be mathematically stated as:
\begin{subequations}
\begin{align}
&\underset{K,\tau_d,\tau_o}{\text{minimize}} \ \ \  J(K,\tau_d,\tau_o) + g(K), \\
 &\text{subject to} \ \ \  K \in \mathbb{K},\\
 & \hspace{1.75cm} S(\tau_d,\tau_o,K) \leq S_B,
  \end{align}
\end{subequations}
where $S$ is given by \eqref{bw2}, and $g(K)$ is a sparsity-promoting function which will be introduced in Section \ref{subsection:innerloop}. The closed-form expression of $J$ is derived next.
\subsection{ \texorpdfstring{$\mc{H}_2$}{H2} norm for the Delayed System}
\noindent The delayed system \eqref{delayed} is infinite dimensional. In order to obtain a linear, finite dimensional LTI approximation of \eqref{delayed}, we use the method of spectral discretization given in \cite{spectral}. Since $\tau_o > \tau_d$ in \eqref{delayed}, following \cite{spectral}, we divide $[-\tau_o,0]$ into a grid of $N$ scaled and shifted Chebyshev extremal points
\begin{equation}
\theta_{k+1} = \frac{\tau_o}{2}\left(\cos\left(\frac{(N-k-1)\pi}{N-1}\right) -1\right), \ k=\{0,\ldots,N-1\}, \label{theta}
\end{equation}
such that $\theta_1=-\tau_o$ and $\theta_N=0$. {The choice of $N$ is guided by \cite[Section 4]{spectral}}. Let $\upsilon(\theta)=x(t+\theta)$ denote the $\theta$-shifted state vector. The extended state  $\eta$ and the closed-loop state matrix $A_{cl}$ can then be written as:
\begin{subequations}
\label{spectral}
\begin{align}
&\eta = [\upsilon^T(\theta_1),  \cdots ,\upsilon^T(\theta_N)= x(t)]^T, \ l_j(\theta)=\prod\limits_{m=1,\ m\neq j}^N \frac{\theta -\theta_m}{\theta_j -\theta_m}, \label{lj}\\
& {A_{cl}}_{ij} = \begin{cases}
 \partial_\theta l_{j}(\theta_i) I_n, \ \ \ j=1,\ldots,N, \ i=1,\ldots,N-1\\
 l_N(-\tau_d) B K_d +  A, \ \ \ \  j=N, \ i=N\\
 l_1(-\tau_d)B K_d + B  K_o, \ \ \ \  j=1, \ i=N\\
 l_j(-\tau_d) B K_d , \ \ \ \ \ j=2,\ldots,N-1, \ i=N,
\end{cases} \label{atilde}
\end{align}
\end{subequations}
where $K_d = K\circ \mc{I}_d$, $K_o = K \circ \mc{I}_o$. 
We can separate $A_{cl}$ into three sub-components:
\begin{align}
&\hspace{1.5cm} A_{cl} = \tilde{A} - \mc{B} K_o N^T_o - \mc{B} K_d N^T_d, \label{acl}\\
&\mc{B}=M B, \ M=[\bb{0},\ldots,\bb{0},I_n]^T, \ N_o = [I_n,\bb{0},\ldots,\bb{0}]^T,
\end{align}
where the first sub-component $\tilde{A}$ is independent of $K_d$ and $K_o$, the second is only dependent on $K_o$, and the third on $K_d$. The explicit expressions for $\tilde{A}$ and $N_d$ in terms of $\tau_d$ and $\tau_o$ will be derived shortly in the next section. The linear approximation of the closed-loop system \eqref{delayed} becomes:
\begin{subequations}
\label{extendeddelayed}
\begin{align}
&\dot{\eta}(t) = A_{cl} \eta(t) + \mc{B}_w w(t),  \\
&  z(t) = \mc{C}\eta(t), \ \mc{C} = \begin{bmatrix}
Q^{\nicefrac{1}{2}} M^T \\ -R^{\nicefrac{1}{2}}(K_d N^T_d + K_o N^T_o)
\end{bmatrix},
\end{align}
\end{subequations}
where  $\mc{B}_w = M B_w$. The algebraic Riccati equations (AREs) and the closed-loop $\mc{H}_2$-norm $J$ can be written as:
\begin{align}
&A_{cl}^T P + P A_{cl} = -\mc{C}^T\mc{C} =-\big( \tilde{Q} + \tilde{C}^T R \tilde{C}\big), \label{CC}  \\
&A_{cl} L + L A_{cl}^T = -\mc{B}\mc{B}^T ,  \label{LL}\\
&J(K,\tau_d,\tau_o) = \text{Tr} ( \mc{B}^T P \mc{B} ) = \text{Tr} (\mc{C} L \mc{C}^T). \label{J}
\end{align}
where $\tilde{Q} = M Q M^T$ and $\tilde{C}= (K_d N^T_d + K_o N^T_o)$.

\section{ Derivation of the gradient of \texorpdfstring{$\mc{H}_2$}{H2} norm}
\label{sec:proposedcodesign}

\noindent Our goal is to design $(K,\,\tau_d,\,\tau_o)$ to minimize $J$. However, from (\ref{CC})-(\ref{J}),  we see that $J$ is a function of $\tilde{A}$ and $N_d$, besides $K$. To compute the gradient of $J$ with respect to $(K,\,\tau_d,\,\tau_o)$, it is important to express $\tilde{A}$ and $N_d$ in terms of these three design variables. We begin this section with these derivations as follows.

\subsection{\texorpdfstring{$\mc{H}_2$}{H2} Performance and Design Variables}

\noindent Recall that the closed-loop state matrix $A_{cl} = \tilde{A} - \mc{B} (K_o N^T_o$ $ +K_d N^T_d  )$. In the next two lemmas, we express $A_{cl}$ as a function of $\tau_o$, $K$ and the {delay} ratio $c=\nicefrac{\tau_d}{\tau_o}$.
\begin{lemma}
\label{lemmatau}
$\tilde{A}$ is a function of $\tau_o$, and can be written as:
\begin{align}
\tilde{A} = \frac{1}{\tau_o}\Lambda + \bar{A}, \ \bar{A}=\texttt{Diag}(\bb{0},A),
\end{align}
where $\Lambda$ is a constant matrix for constant $N$. $\sqbullet$
\end{lemma}
\begin{lemma}
\label{lemmand}
$N_d$ is a function of the ratio $c=\nicefrac{\tau_d}{\tau_o} \in [0,1]$, and can be written as:
\begin{equation}
N_d (c) = \left(\Gamma \nu (c)\right) \otimes I_n, \ \nu (c) = [c^{N-1} \ c^{N-2} \ \ldots \ c^2 \ c \ 1]^T, \label{Nddeclare}
\end{equation}
where $\Gamma \in \mathbb{R}^{N\times N}$ is a constant matrix for constant $N$. $\sqbullet$
\end{lemma}
\smallskip
\par Lemmas \ref{lemmatau} and \ref{lemmand} show that {for fixed $N$}, $J$ for the system in \eqref{extendeddelayed} can be written as a function of $\tau_o$ and $c$. Henceforth, all of our analysis for minimizing $J$ will be carried out using $\tau_o$ and $c$, instead of $\tau_o$ and $\tau_d$. This change of variables is invertible, and therefore, there is no loss of generality.

\subsection{Gradient of \texorpdfstring{$\mc{H}_2$}{H2} norm}
In order to minimize $J$, we next derive the gradient of $J$. We define a set $\mc{K}$ as: 
\begin{align}
 \mc{K}:= \{(K,\tau_o,c): \te{Re}\big(\lambda(A_{cl})\big) < 0 \},   \end{align}
i.e., the set of solutions that guarantee closed-loop stability of \eqref{extendeddelayed}. Given this definition, we first prove the existence of a unique solution of \eqref{CC} and differentiability of $P$, followed by the derivation of $\nabla J$. For the rest of the paper, the $A^{'}(B)$ notation represents differentiability of $A$ depending on $B$.
\begin{lemma}
\label{lemmaunique}
{ Let $(K,\tau_o,c)\in \mc{K}$. Then, there exists a unique solution $P(K,\tau_o,c)$ of \eqref{CC}. Moreover, $P$ is differentiable with respect to the variables $\tau_o$, $c$ and $K$ on $\mc{K}$. Specifically, $P^{'}(\tau_o) d\tau_o$, $P^{'}(c) dc$ and $P^{'}(K) dK$ follow as solutions of the following Lyapunov equations, respectively:}
\begin{align}
&A_{cl}^T P^{'}(\tau_o) d\tau_o + P^{'}(\tau_o) d\tau_o A_{cl} = \frac{d\tau_o}{\tau^2_o} (\Lambda^T P + P\Lambda ) ,\label{taudiff}\\
& A_{cl}^T P^{'}(c) dc + P^{'}(c) dcA_{cl} =  N^{'}_d dc K^T_d G + G^T  K_d {N^{T}_d}^{'}dc , \label{cdiff}\\
&A_{cl}^T P^{'}(K) dK + P^{'}(K) dKA_{cl} =  -Z_d - Z^T_d - Z_o - Z^T_o,  \label{psiK}
 \end{align}
where $G=\big( R (K_d N^T_d(c)$ $+ K_o N^T_o) - \mc{B}^T P \big)$, $Z_d =N^T_d (dK$ $ \circ \mc{I}_d)G$, $Z_o =N^T_o (dK \circ \mc{I}_o) G$ and $N^{'}_d = \left(\Gamma \partial \nu(c)\right)\otimes I_n$. $\sqbullet$
\end{lemma}
\smallskip
\par { We next use Lemma \ref{lemmaunique} to state the following theorem.}
\smallskip
\begin{theorem}
\label{Theoremgrad}
 $J$ in \eqref{J} is differentiable on $\mc{K}$. The gradient of $J$ is evaluated as:
\begin{align}
&J'(\tau_o)  = -\frac{2}{\tau^2_o}\te{Tr} ( \Lambda^T P L), \  J'(c) =  2 \te{Tr}( N^{'}_d K^T_d G L), \label{JtauJc} \\
&\nabla J (K) = 2 ( (G L N_d)\circ \mc{I}_d + (G L N_o)\circ \mc{I}_o). \label{J(K)}
\end{align}
% where $L$ is the solution of \eqref{LL}.
\end{theorem}
The negative directions of $ J'(c)$ and $J'(\tau_o) $, as derived in Theorem \ref{Theoremgrad}, always point to the trivial solution $c=0, \ \tau_o=0$ which defeats the purpose of designing $\tau_d$ and $\tau_o$. This is because the partial derivatives in \eqref{JtauJc}-\eqref{J(K)} are derived with the assumption that $K$, $\tau_o$ and $c$ are independent of each other as $K^{'}(\tau_o)$ and $K^{'}(c)$ cannot be computed directly given the implicit dependence of $K$ on $\tau_o$ and $c$. Therefore, it would be incorrect to co-design $c$, $\tau_o$ and $K$ using just the gradient information. Starting from a stabilizing $(K,\tau_o,c)\in\mc{K}$, as soon as we change either $\tau_o$ or $c$, we must update $K$ to ensure stability of \eqref{extendeddelayed}. In other words, $(K,\tau_o)$ and $(K,c)$ must be co-designed separately in sequence while holding $c$ and $\tau_o$ as constant in the respective steps.
\subsection{Co-design of Controller and Delays}
 We next describe how equations in \eqref{CC}-\eqref{LL} can be relaxed for each of the two co-designs.
\par $\bullet$ \textit{\textbf{Co-design of $(K,\tau_o)$}}
\begin{theorem}
\label{theoremktau}
Let $\omega_o= \nicefrac{1}{\tau_o}$. Consider a known tuple $(K^*,\omega^*_o,c^*)\in\mc{K}$ satisfying \eqref{CC} with a known $P^*$ for closed-loop state matrix $A^{*}_{cl}(K^{*},\omega^{*},c^{*})$. Let ${\omega_o} = \omega^*_o + {\Delta \omega}$, ${K}=K^* + {\Delta K}$, ${P}=P^* + {\Delta P}$ and ${\alpha}\in\mathbb{R}$ be obtained as a solution of the following SDP.
% is a stabilizing solution for \eqref{extendeddelayed}:
\begin{subequations}
\label{theoremktaueq}
\begin{align}
&{\phi_0} + {\phi_1} + {\psi_0} + {\alpha} I \succeq 0, \\
& |{\Delta \omega}| \leq \zeta_1,\ \|{\Delta P} \| \leq \zeta_2, \label{39bc} \\
& {\alpha} \geq 2\zeta_1 \|\Lambda^T {\Delta P}\| + 2\zeta_2 \|\mc{B}{\Delta \tilde{C}}\| + \| R^{\nicefrac{1}{2}}{\Delta \tilde{C}}\|^2, \label{39d}
\end{align}
\end{subequations}
where $\alpha$, $\Delta K$, $\Delta P$ and $\Delta \omega$ are the design variables. Then, $(K,\nicefrac{1}{\omega_o},c^{*})$ is a stabilizing tuple for \eqref{extendeddelayed}. In \eqref{theoremktaueq}, ${\phi_0} = A_{cl}^{*T}{ P} + {P}A^*_{cl}$, $K^*_d = K^* \circ \mc{I}_d$, ${\Delta K_d} = {\Delta K} \circ \mc{I}_d$, $K^*_o = K^* \circ \mc{I}_o$, ${\Delta K_o }= {\Delta K} \circ \mc{I}_o$, $\tilde{C}^* = \tilde{Q} + (K^*_d N^T_d + K^*_o N_o^T)$, ${\Delta \tilde{C}}=({\Delta K_d} N^T_d +{\Delta K_o} N^T_o)$, ${A_1 }= -\mc{B} ({\Delta \tilde{C}}) + {\Delta \omega} \Lambda$, ${\phi_1 }= {A}^T_{{1}} P^* + P^* {A_1}$, ${\psi_0} = \tilde{C}^{*T} R \tilde{C}^* + {\Delta \tilde{C}}^T R \tilde{C}^* + \tilde{C}^* R {\Delta \tilde{C}}$  and,  $\zeta_1, \ \zeta_2$ are chosen constants. $\sqbullet$
\end{theorem}

\par $\bullet$ \textit{\textbf{Co-design of $(K,c)$}}
\par Next, consider the co-design step for $(K,c)$. Recall that $A_{cl}$ is a non-linear function of $c \in [0,1]$ through $N_d(c)$ as shown in Lemma \ref{lemmand}, and therefore, the exact expression of $N_d(c)$ cannot be used while forming the SDP relaxations. To circumvent this problem, we divide $[0,1]$ into $k_c$ sub-intervals $[{c_1},{c_2}],[{c_2},{c_3}],\ldots,[c_{k_c},c_{k_c + 1}]$ with each sub-interval small enough to allow $N_d(c)$ to be approximated as an affine function $\hat{N}_d(c)$. Let each sub-interval $[c_i,c_{i+1}]$ have an associated $\chi^{(i)} \in \mathbb{R}^{N\times 2}$ as the vector of affine coefficients. The approximated function is written as:
\begin{align}
&\hat{N}_d (c) = \left(\chi^{(i)} [c, 1]^T\right)\otimes I_n, \ c\in [c_i,c_{i+1}], \ i=1,\ldots,k_c.
\end{align}
The coefficients can be computed from a linear curve fitting on \eqref{Nddeclare}. Larger the number of sub-intervals $k_c$, lower is the approximation error $\|\hat{N}_d - N_d\|$. For our simulations in Section \ref{sec:simulations}, we have used $k_c=10$. We next present the SDP relaxation for the co-design of $(K,c)$.
\begin{theorem}
\label{theoremkc}
Consider a known tuple $(K^{*},\tau^{*}_o$, $c^{*})\in\mc{K}$ with $c^*\in[c_i,c_{i+1}]$ for some $i\in \{ 1,\ldots,k_c\}$ satisfying \eqref{LL} with a known $L^{*}$ for closed-loop state matrix $A^{*}_{cl}(K^{*},\tau^{*}_o,c^{*})$. Let $c=c^{*} + \Delta c$, $K=K^{*} + \Delta K$, $L=L^{*}$ $+\Delta L$ and $\alpha\in\mathbb{R}$ be a solution of the following SDP:
\begin{subequations}
\label{theoremkceq} 
\begin{align}
&\phi_0 + \phi_1 + \mc{B}\mc{B}^T + \alpha I \succeq 0, \\
& c_i \leq c\leq c_{i+1}, \ \|\Delta L\| \leq \beta, \label{beta}\\
& \alpha \geq 2 \beta \|\mc{B} (\Delta K_d N^T_d (c^{*}) + \Delta K_o N^T_o)\| + 2\beta \mathfrak{S} \|\mc{B}\Delta K_d \| \nn \\
& + 2 \beta \|\Delta N_d K^{* T}_d B^T M^T\| + 2 \mathfrak{S} \|\mc{B}\Delta K_d\|\|L^*\| , \label{alphaconstraint}
\end{align}
\end{subequations}
where $\alpha$, $\Delta K$, $\Delta P$ and $\Delta c$ are the design variables. Then, $(K,\tau^{*}_o,c)$ is a stabilizing tuple for \eqref{extendeddelayed}. In \eqref{theoremkceq}, $\Delta N_d = \hat{N}_d (c) - N_d (c^*)$, $\phi_0 = A_{cl}^{*} L + L A_{cl}^{* T}$, $\phi_1 = A_1 L^{*} + L^{*} A^T_1$, $A_1 = -\mc{B} (K^{*}_d \Delta N^T_d(c) + \Delta K_d N^T_d (c^{*}) + \Delta K_o N^{T}_o)$. The scalar $\beta$ is a chosen constant, and $\mathfrak{S}\geq \|N_d(c)\|$. $\sqbullet$
\end{theorem}

\par Starting from a known stabilizing tuple $(K^*,\tau^*,c^*)$, Theorems \ref{theoremktau} and \ref{theoremkc} enable us to co-design new stabilizing pairs $(K,\tau_o)$ and $(K,c)$, respectively. Next, we integrate the bandwidth cost constraint \eqref{bw2} with the SDPs in \eqref{theoremktaueq} and \eqref{theoremkceq}.

\subsection{Incorporating Bandwidth Constraints}
\label{subsec:changeofvariables}
\par \noindent We impose the bandwidth cost constraint \eqref{bw2} as part of \textbf{P1}, which can be rewritten as:
\begin{align}
S = 2m_{cp}\left( \frac{N_{row}(K)+N_{col}(K) }{c \tau_o}\right) + m_{cc}\left( \frac{N_{off}(K)}{\tau_o - c\tau_o}\right) \leq S_{b}. \label{constraintchange}
\end{align}
Recall that $S$ is the total bandwidth cost and $S_b$ is the upper bound imposed on it. When \eqref{constraintchange} is imposed on SDPs \eqref{theoremktaueq} and \eqref{theoremkceq}, we obtain an alternative form of \eqref{constraintchange}, which is stated in the next proposition.
\begin{prop}
\label{propositionkcktau}
Consider a known tuple $(K^{*},\tau^*_o,c^*)\in\mc{K}$ with an associated bandwidth cost $S^*\leq S_b$. Denoting $n^*_{cp}=N_{row}(K^*)+N_{col}(K^*)$ and $n^*_{cc}=N_{off}(K^*)$, the following statements are true.\\
1) Keeping $\tau_o=\tau^*_o$, let $c^*$ be perturbed to $c$ resulting in a cost $S$. Then, $\delta S(c): = S-S^*$ is a convex function of $c$:
\begin{align}
 \delta S(c) = \frac{(S^{*} \tau^*_o) c^2 + (m_{cc} n^*_{cc}- 2m_{cp}n^*_{cp} - S^{*} \tau^*_o )c + 2 m_{cp}n^{*}_{cp}}{c(1-c)\tau^{*}_o}. \label{constraintchange2}  
\end{align} 
The constraint $\delta S(c) \leq 0$ implies $S\leq S_b$.\\
2) Keeping $c=c^*$, let $\tau^{*}_o$ be perturbed to $\tau_o$, resulting in a new bandwidth cost $S$. Then, $\delta S (\tau_o):= S-S^*$ is an affine function of $\tau_o$:
\begin{align}
\delta S(\tau_o) = \frac{1}{S^*}\left(\frac{2 m_{cp} n^*_{cp}}{c^*}  + \frac{m_{cc} n^*_{cc}}{(1-c^*)}\right) -\tau_o. \label{constraintchange3}
\end{align}
The constraint $\delta S(\tau_o)\leq 0$ implies $S\leq S_b$.
\end{prop}
\par \noindent \textit{Proof:} The proof follows from simple algebra. $\sqbullet$
\par Since $\delta S (\tau_o)$ and $\delta S(c)$ are each convex in their respective arguments in the above proposition, we can easily incorporate them in the co-design SDPs of Theorems \ref{theoremktau} and \ref{theoremkc} to satisfy the bandwidth constraint in \eqref{constraintchange}. Note that since $K$ is co-designed with either $\tau_o$ or $c$, the true bandwidth cost $S$ depends on $K$ as well through $(N_{row}(K)+ N_{col}(K))$ and $N_{off}(K)$. If $N_{row}(K)\leq N_{row}(K^*)$, $N_{col}(K)$ $\leq N_{col}(K^*)$ and $N_{off}(K)\leq N_{off}(K^*)$, one can easily verify that $\delta S(c)\leq 0$ and $\delta S(\tau_o)\leq 0$ in \eqref{constraintchange2}-\eqref{constraintchange3} hold, and the true bandwidth costs always satisfy \eqref{constraintchange}. We ensure this fact by imposing a two-loop structure in our design algorithm, as will be seen shortly in the next section. We next bring together the co-design SDPs \eqref{theoremktaueq}, \eqref{theoremkceq} and bandwidth constraints \eqref{constraintchange2}, \eqref{constraintchange3} in the form of our main algorithm.

\section{Problem Setup in Two-Loop ADMM Form}
\label{sec:problemsetupinadmm}
The $\mc{H}_2$-norm $J$, in general, increases with increasing sparsity of $K$ \cite{negi2018}, while the bandwidth cost $S$ reduces. Due to these inherent trade-offs between the objectives and the constraints, \textbf{P1} is a prime candidate to be reformulated as a two-loop ADMM optimization. The outer-loop co-designs $(K,\tau_o)$ and $(K,c)$ using \eqref{theoremktaueq}-\eqref{theoremkceq} under the  bandwidth constraints \eqref{constraintchange2}-\eqref{constraintchange3}. The inner-loop, on the other hand, sparsifies $K$ while minimizing $J$. We describe the inner and outer loops in Sections \ref{subsection:innerloop} and \ref{subsection:outerloop} respectively, followed by the main algorithm in Section \ref{subsection:mainalgo}.
\subsection{Inner ADMM Loop}
\label{subsection:innerloop}
\noindent Throughout the inner ADMM loop, we hold both $\tau_o$ and $c$ as constants. The mathematical program of the inner loop denoted as \textbf{P1}\textsubscript{in} is written as follows:
\begin{subequations}
\label{admmform}
\begin{align}
\text{\textbf{P1}\textsubscript{in}} \ : \ \ & \minimize_{K,F} \ \ J(K) + \lambda g(F), \\
&\text{subject to} \ \ \ K=F,
\end{align}
\end{subequations}
where $\lambda$ is a regularization parameter and $g(F)=\|W\circ F\|_{l_1}$ is the weighted $l_1$ norm function which is used to induce sparsity in $F$. The weight matrix $W$ for $g(F)$ is updated iteratively through a series of reweighting steps from the solution of the previous iteration as \cite{weightedl1}:
\begin{equation}
    W_{ij} = \frac{1}{|F_{ij}| +\epsilon }, \ \ 0<\epsilon \ll 1. \label{wl1}
\end{equation}
The augmented Lagrangian for \textbf{P1}\textsubscript{in} is 
\begin{equation}
\mc{L}_p = J(K) + \lambda g(F) + \text{Tr}(\Theta^T (K-F)) + \frac{\rho}{2}\|K-F\|^2_{\te{F}},\label{auglag}
\end{equation}
where $\rho$ is a positive scalar and $\Theta$ is the dual variable. ADMM involves solving each objective separately while simultaneously projecting onto the solution set of the other. As shown in \cite{mihailo,boyd}, \eqref{auglag} is used to derive a sequence of iterative steps $K$-min, $F$-min and $\Theta$-min by completing the squares with respect to each variable.
\begin{subequations}
\label{admmeq}
\begin{align}
& K_{k+1} = \underset{K}{\te{argmin}} \ \Phi_1(K)=\underset{K}{\te{argmin}} \ J(K) + \frac{\rho}{2} \|K-U_k\|^2_{\te{F}}, \label{Kmin} \\
& F_{k+1}= \underset{F}{\te{argmin}} \ \Phi_2(F) =\underset{F}{\te{argmin}} \ \lambda g(F) + \frac{\rho}{2} \|F-V_k\|^2_{\te{F}},\label{Fmin}\\
& \Theta_{k+1} = \Theta_k + \rho(K_{k+1} - F_{k+1}),\label{Theta}
\end{align}
\end{subequations}
where $U_k = F_k - \frac{1}{\rho} \Theta_k $ and $V_k = K_{k+1}+\frac{1}{\rho}\Theta_k$. We next present a method to solve $K$-min and provide an analytical expression for $F$-min.
\subsubsection{\texorpdfstring{$K$}{K}-min Step}
Setting $\nabla \Phi_1(K)=0$ and using Theorem \ref{Theoremgrad}, we get the following condition for optimality\footnote{We denote $N_d(c)$ simply as $N_d$ throughout this subsection as $c$ is constant for \textbf{P1}\tsb{in}.}:
\begin{equation}
\left[ (G L N_d)\circ \mc{I}_d+ (G L N_o)\circ \mc{I}_o \right] + \frac{\rho}{2}(K - U) =0, \label{kmineqgrad}
\end{equation}
where $G = R(K_d N^T_d + K_o N^T_o) -\mc{B}^T P$ and $U=U_k$ for the $(k+1)$-th iteration of the ADMM loop. $P$ and $L$ are the solutions of AREs \eqref{CC} and \eqref{LL}, respectively. $K$-min begins with a stabilizing $K$, solves \eqref{CC} and \eqref{LL} for $P$ and $L$, and then solves \eqref{kmineqgrad} to obtain a new gain $\bar{K}$ as follows: 
\begin{align}
&\bar{K} = \texttt{Reshape}\left( (\hat{V}_d \circ T_d  +  \hat{V}_o \circ T_o  + \rho I_{n^2})^{-1} \mu   ,\texttt{[m,n]}\right),\label{tddactual}\\
&T_d = (T_{dd} \circ \hat{V}^T_d+ T_{od}\circ \hat{V}^T_o), \  T_o = (T_{oo}\circ \hat{V}^T_o + T_{do} \circ \hat{V}^T_d ),  \nn\\
&T_{dd} = 2 (N^T_d L N_d \otimes R), \  T_{od} = 2(N^T_o L N_d \otimes R), \nn \\
&T_{oo} = 2 (N^T_o  L N_o \otimes R ), \ T_{do} = 2(N^T_d L N_o \otimes R), \nn \\
& \mu = \texttt{vec} \left( (2 \mc{B}^T P L N_d) \circ \mc{I}_d +  (2 \mc{B}^T P L N_o) \circ \mc{I}_o + \rho U \right),\nn \\
& \hat{V}_d = \mathbf{1} \otimes v_d, \ v_d = \texttt{vec}(\mc{I}_d),\ \hat{V}_o = \mathbf{1} \otimes v_o, \ v_o = \texttt{vec}(\mc{I}_o), \nn
\end{align}
The notation $B=\texttt{Reshape}(A,[p,q])$ is used for an opera-tor that reshapes $A\in\mathbb{R}^{m\times n}$ in row-traversing order to another matrix $B\in\mathbb{R}^{p\times q}$, provided $pq=mn$. We use $\texttt{vec}$ to represent the vectorization operator and $\mathbf{1}^T\in\mathbb{R}^{n^2}$ to represent a vector of all ones. For details of the above derivation, see the Appendix. It can be shown that $\tilde{K} = K-\bar{K}$ is the descent direction for $\Phi_1$ \cite[See Lemma 4.1]{computational}. The Armijo-Goldstein line search method can then be used to determine a step size $s$ to ensure $(K+s\tilde{K}) \in \mathbb{K}$, { i.e., stability of \eqref{extendeddelayed} is maintained}. The iterative process continues till we obtain $\nabla \Phi_1(K)\approx 0$.

\subsubsection{\texorpdfstring{$F$}{F}-min Step}
\label{subsub1}
The solution of the $F$-min step is well-known in the literature \cite[Sec. 4.4.3]{boyd} as:
\begin{equation}
F_{ij} = \begin{cases}
(1- \frac{a_{ij}}{|V_{ij}|} )V_{ij}, \ \te{if} \ |V_{ij}| > a_{ij}, \\
0, \ \te{otherwise},
\end{cases}
\end{equation}
where $a_{ij}=\frac{\lambda}{\rho} W_{ij}$. Note that large values of $\lambda$ will induce more sparsity, and therefore may lead to a sudden increase in $J$. Therefore, $\lambda$ must be increased in small steps. { The regularization path, for example, can be logarithmically spaced from $0.01\lambda_{max}$ to $0.95\lambda_{max}$, where $\lambda_{max}$ is ideally the critical value of $\lambda$ above which the solution of \textbf{P1}\textsubscript{in} is $K=F=0$ \cite{boyd}. In our simulations, $\lambda_{max}=1$.}
\subsection{Outer Loop}
\label{subsection:outerloop}
\noindent The outer-loop of our algorithm designs $\tau_o$ and $c$ with bandwidth constraint \eqref{constraintchange} and updates the weight matrix $W$ for minimizing the weighed $l_1$ norm in \eqref{wl1}. Co-design of $K$ in this loop is necessary to ensure stability as $\tau_o$ and $c$ change. Let $K^*=F^*$ and $\Theta^*$ be the output of the last converged inner loop with $U^* = K^* - \frac{1}{\rho}\Theta^* $. Programs \textbf{P1}\tsb{o1} and \textbf{P1}\tsb{o2} directly design $(K,\tau_o)$ and $(K,c)$, respectively, in sequence as follows: 
\begin{subequations} 
\label{p101p102}
\begin{align}
\te{\textbf{P1}\tsb{o1}}\ : \  &\underset{K,\tau_o,P}{\te{minimize}} \ \hat{J}(K,\tau_o) + \frac{\rho}{2} \|K-U^*\|^2_{\te{F}},\\
 &\te{subject to } \ \delta S(\tau_o) \leq 0, \\ &\hspace{1.65cm} \text{SDP in }\te{Eq.} \  \eqref{theoremktaueq},  \\ 
\te{\textbf{P1}\tsb{o2}} \ :  \ &\underset{K,c,L}{\te{minimize}} \ \hat{J}(K,c) + \frac{\rho}{2} \|K-U^*\|^2_{\te{F}},\\
 &\te{subject to } \ \delta S(c) \leq 0, \\
 &\hspace{1.65cm} \text{SDP in }\te{Eq.} \  \eqref{theoremkceq}, 
\end{align}
\end{subequations}
where $\hat{J}(K,\tau_o)=\te{Tr} (\mc{B}^T P \mc{B})$, $\hat{J}(K,c) $ $=\te{Tr}(L \mc{C}^{*T}\mc{C}^* )$, $\mc{C}^*$ $=( K^*\circ \mc{I}_d) N^T_d(c^*) + (K^*\circ\mc{I}_o)N^T_o$. We next present our main algorithm to show the iterative solutions of \textbf{P1}\tsb{o1} and \textbf{P1}\tsb{o2} beginning from a known stabilizing tuple $(K^*,\tau^*_o,c^*)$.
\begin{algorithm}[t]
\caption{Main Algorithm}
\label{MainAlgorithm}
\SetAlgoLined
 \nl\KwIn{Initial feasible point $(K^*_{o},\tau^*_o$, $c^*)\in\mc{K}$ }
 \nl\For{$\lambda_i = 0.01\lambda_{max}$ to $0.95\lambda_{max}$}{
  \nl\KwIn{$K^*$, $\tau^*_o$ and $c^*$ stabilizing for \eqref{delayed}  }
  \nl\For{$1$ to Maximum Reweighted Steps}{
   \nl Solve \textbf{P1}\tsb{o1} using $K^*$, $\tau^*_o$, $c^*$ to get $\hat{K}$, $\tau_o$ \\
  \nl Solve \textbf{P1}\tsb{o2} using $\hat{K}$, $\tau_o$, $c^*$ to get updated $K$, $c$.\\
 \nl\KwIn{Inner loop initial: $K$, $c$, $\tau_o$ }   
 \nl \label{terminating}\While{ADMM Stopping Criteria \textbf{not} met} 
 { \nl\noindent $\text{$K${-min : }}$ Solve \eqref{Kmin} for $K_{k+1}$\\
 \nl\noindent $ \text{$F${-min : }}$ Solve \eqref{Fmin} for $F_{k+1}$\\ 
\nl\noindent Update $\Theta$ using \eqref{Theta}
}
\nl\KwResult{ $K^*=K$, $\tau^*_o = \tau_o$, $c^* = c$.}
\nl\noindent Update $W$ using $K^*$ from \eqref{wl1}.\\
\nl\noindent Update $S^*$ using $N_{row}(K^*)$, $N_{col}(K^*)$, $N_{off}(K^*)$, $\tau^*_o$ and $c^*$ from \eqref{constraintchange}.
}
}
\nl\KwResult{$K$, $\tau_o$ and $c$ are obtained for $\lambda_i$.}
\end{algorithm}
\vspace{-0.4cm}
\subsection{Main Algorithm}
\label{subsection:mainalgo}
\noindent Our main algorithm is listed in Algorithm \ref{MainAlgorithm}; the following points explain its key steps.
\smallskip
\par\noindent$\bullet$ Using \textbf{P1}\tsb{o1}, we first co-design a stabilizing pair $(\hat{K},\tau_o)$ from an initial tuple $(K^*, \tau_o^*, c^*)\in\mc{K}$. The two are designed together as the initial $K^*$ may not be stabilizing for $\tau_o$ satisfying the bandwidth constraint \eqref{constraintchange3}.
\par \noindent $\bullet$ We then use the solution of \textbf{P1}\tsb{o1}, i.e., $(\hat{K},\tau_o,c^*)\in\mc{K}$ as the initial point for \textbf{P1}\tsb{o2} to find an updated pair $(K,c)$. From Proposition \ref{propositionkcktau}, $\delta S(c)$ in \eqref{constraintchange3} is convex in $c$. Let $c_{min}$ be the minimizer of $\delta S(c)$. If $\hat{K}$ is stabilizing for $c_{min}$, then instead of co-designing $(K,c)$, we can directly set $c=c_{min}$ and $K=\hat{K}$, and then use a procedure similar to $K$-min to minimize $J(K)$ starting from $\hat{K}$. 
\par \noindent $\bullet$ The inner-loop begins with $(K,\tau_o,c)\in\mc{K}$. $K$ is updated in the direction of decreasing $J$ and increasing sparsity while { $\tau_o$ and $c$} remain constant.
\par \noindent $\bullet$ Following \cite[Sec. III-D]{mihailo} and \cite[Sec. 3.4.1]{boyd}, $\rho$ in \eqref{admmeq} is chosen to be sufficiently large to ensure the convergence of the inner ADMM loop. Since $J$ is nonconvex, convergence of this loop, in general, is not guaranteed, as is commonly seen in the sparsity promoting literature \cite{mihailo}. However, large values of $\rho$ have been shown to facilitate convergence. We use $\rho=100$ for our simulations. The stopping criterion for the inner loop in Line \ref{terminating} of Algo. 1 follows \cite[Sec. 3.3.1]{boyd}.

\section{Simulations Results}
\label{sec:simulations}

\subsection{Delay-Design With No Bandwidth Constraints}
\label{subsec:firstsim}

\begin{figure}
\centering
\includegraphics[scale=0.57]{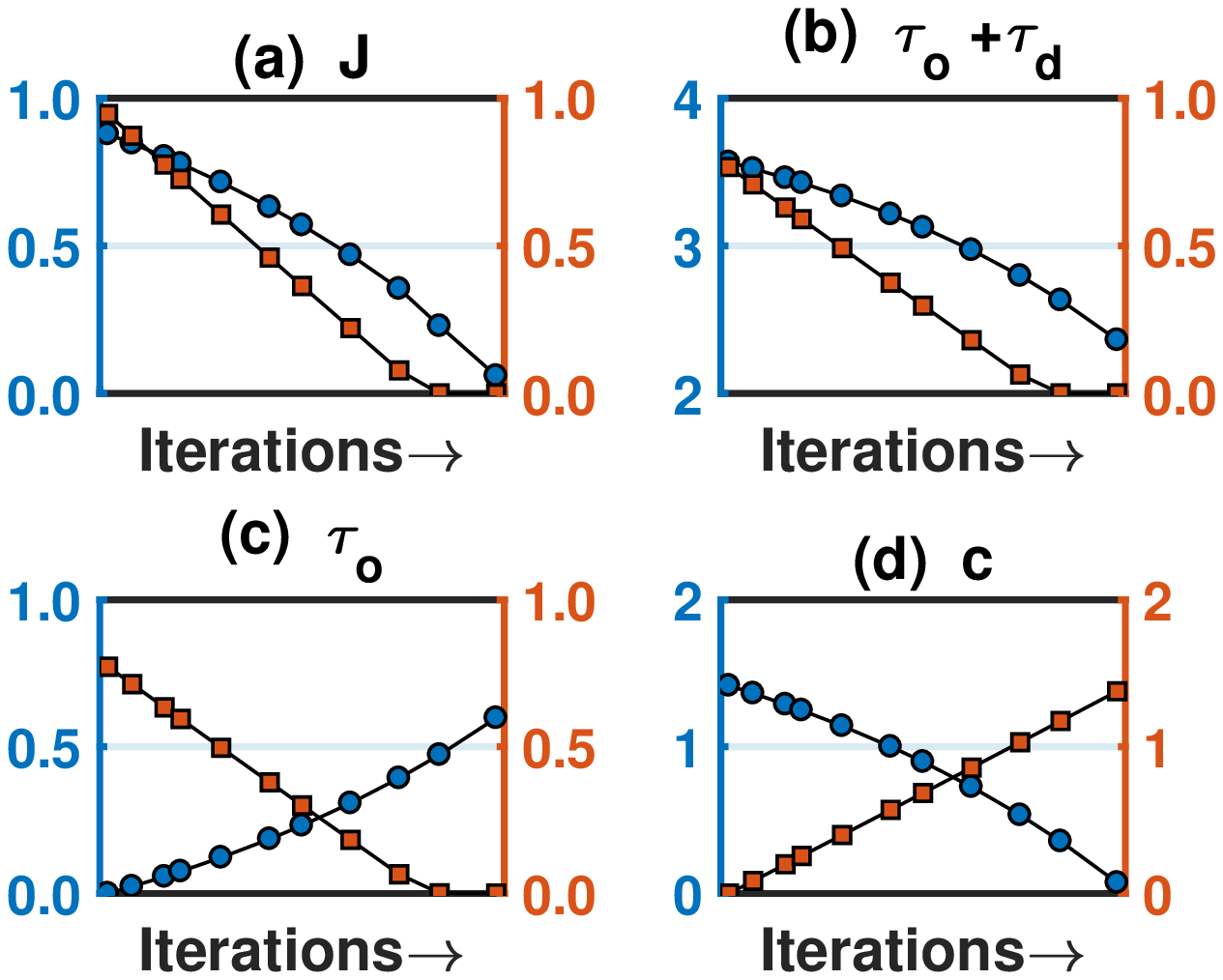}
\caption{(a), (b), (c), (d) show normalized $J$, $\tau_o+\tau_d$, $\tau_o$ and $c$ vs iterations - Right axis for {\color{orange} Model I\tsb{a}} (${\color{orange} \sqbullet}$), Left axis for {\color{cyan}Model I\tsb{b}} (${ \color{cyan}\bullet}$). }
\label{fig:case1}
\end{figure}

\noindent We first present simulations where only the outer loop is iterated without considering any bandwidth constraint in Algorithm \ref{MainAlgorithm}. This example shows that the relative magnitudes of $\tau_d$ and $\tau_o$ for obtaining minimum $\mc{H}_2$-norm can be significantly different for different systems. Absence of the bandwidth cost, as indicated before, will lead to the trivial solution $\tau_o=0$, $\tau_d=0$. To avoid this, we impose a simple artificial constraint $|(\tau_d-\tau^{*}_d) + (\tau_o - \tau^{*}_o)| \leq \epsilon$ where $0<\epsilon\ll 1$ is a small tolerance, and $(K^*,\tau^*_o,\tau^*_d) \in  \mc{K}$ is the initial point for every iteration. This initial tuple is replaced by the newly designed $(K,\tau_o,\tau_d)\in\mc{K}$ at the end of every iteration. We simulate two randomly generated models I\tsb{a} and I\tsb{b} with $A\in\mathbb{R}^{5\times 5}$, $B=B_w=I_n$, $K^*=K_{LQR}$, $Q=R=I_n$ for two different initial conditions as part of Case A. The logarithm of ratios of $J$, $\tau_d+\tau_o$, $\tau_o$ and $c$ with respect to their respective minima are plotted in Fig. \ref{fig:case1}.
\smallskip
\par \noindent \textbf{Case A}: Right and left axis of all the sub-figures in Fig. \ref{fig:case1} show system I\tsb{a} and I\tsb{b} with $(c^*,\tau^*_o)$ chosen as $(0.489,0.141)$ and $(0.833,0.108)$, respectively. For both the systems, $J$ in Fig. \ref{fig:case1} (a) is seen to be decreasing as $\tau_o+\tau_d$ decreases. This is expected as $\mc{H}_2$-performance improves with a decrease in the overall delay. Fig. \ref{fig:case1} (a), (c) and (d) show that for achieving a lower $J$, the model I\tsb{a} requires a lower $\tau_o$ and a higher $c$, while I\tsb{b} requires a higher $\tau_o$ and a lower $c$. We can infer that obtaining a better $\mc{H}_2$-performance can demand completely different relative magnitudes of $\tau_d$ and $\tau_o$ depending on the system model and the initial conditions. Thus, this example validates the motivation of our problem in determining the trade-off between $\tau_d$ and $\tau_o$.
\begin{figure}[hbtp]
\begin{minipage}{0.5\linewidth}
\centering
\includegraphics[scale=0.291]{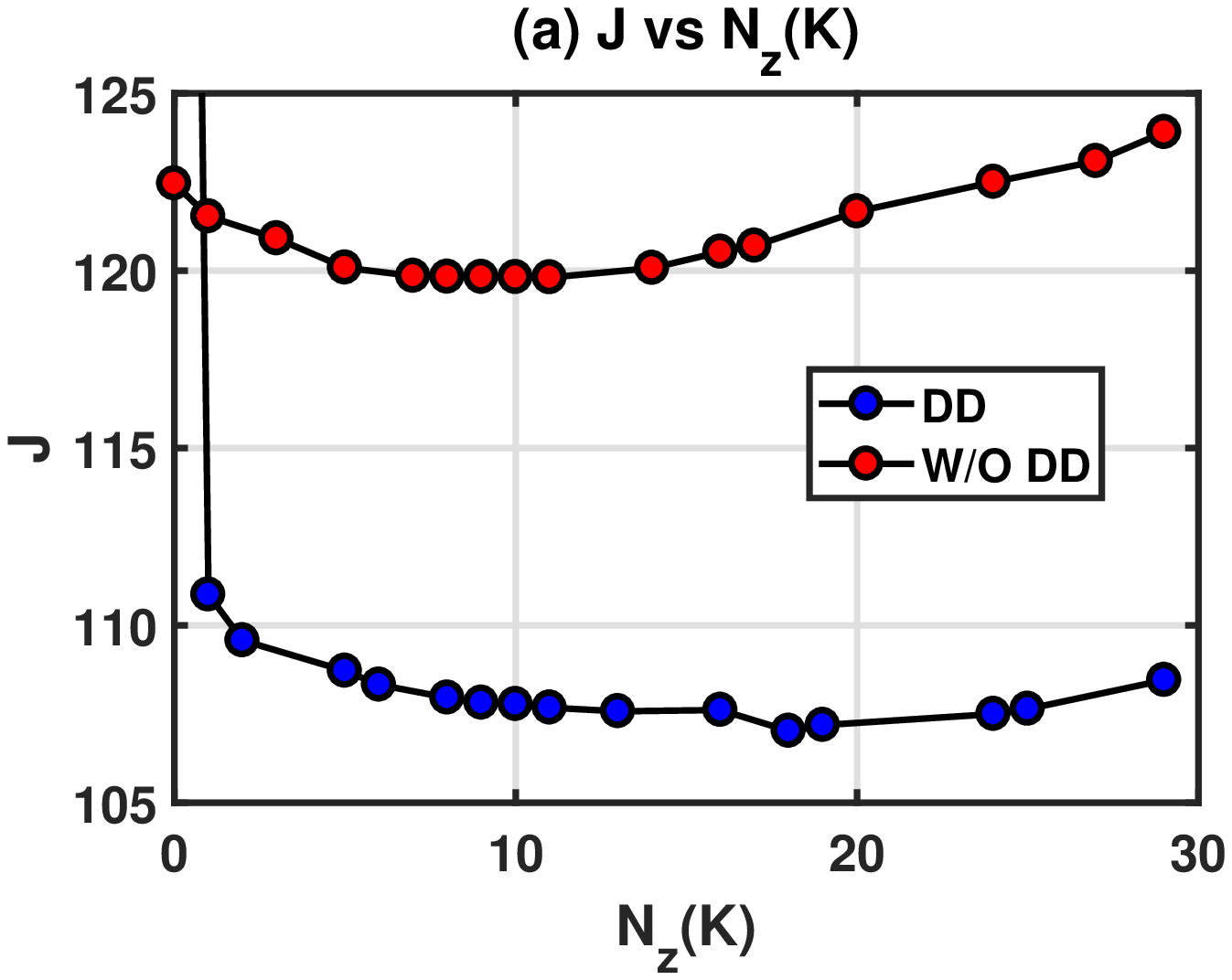} %
\end{minipage}\begin{minipage}{0.5\linewidth}
\centering
 \includegraphics[scale=0.291]{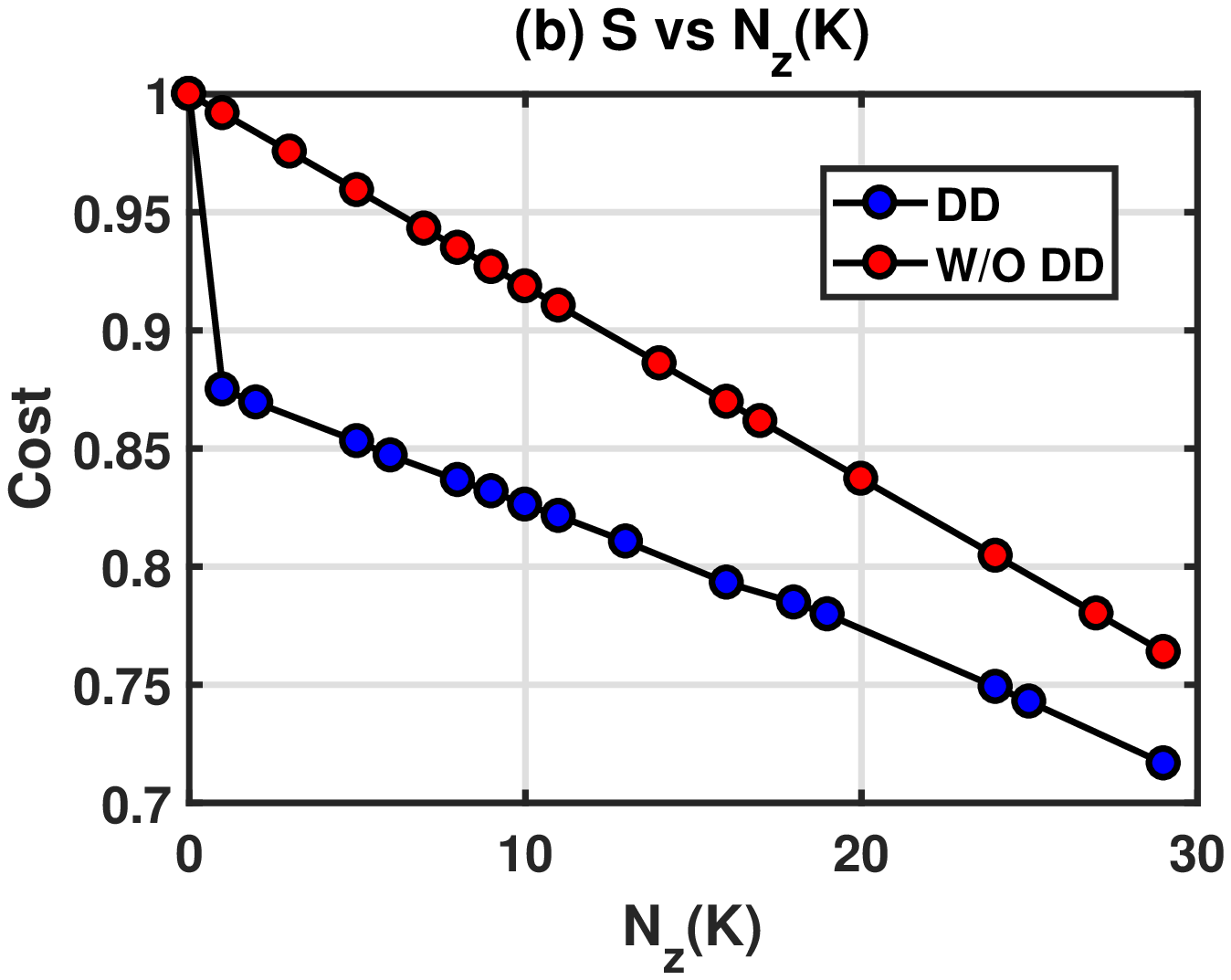}
\end{minipage}\\
\begin{minipage}{0.5\linewidth}
\centering
\includegraphics[scale=0.291]{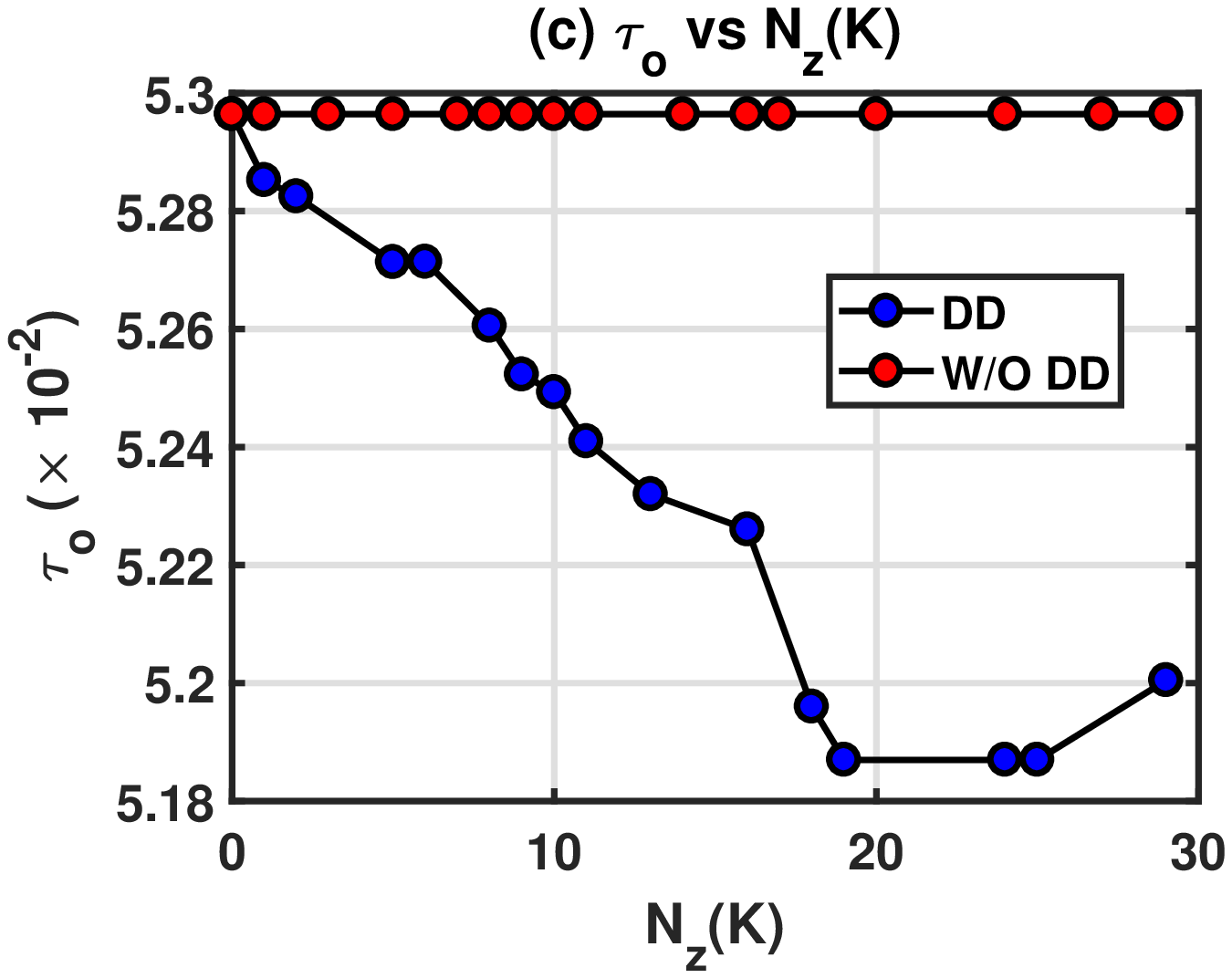} %
\end{minipage}\begin{minipage}{0.5\linewidth}
\centering
 \includegraphics[scale=0.291]{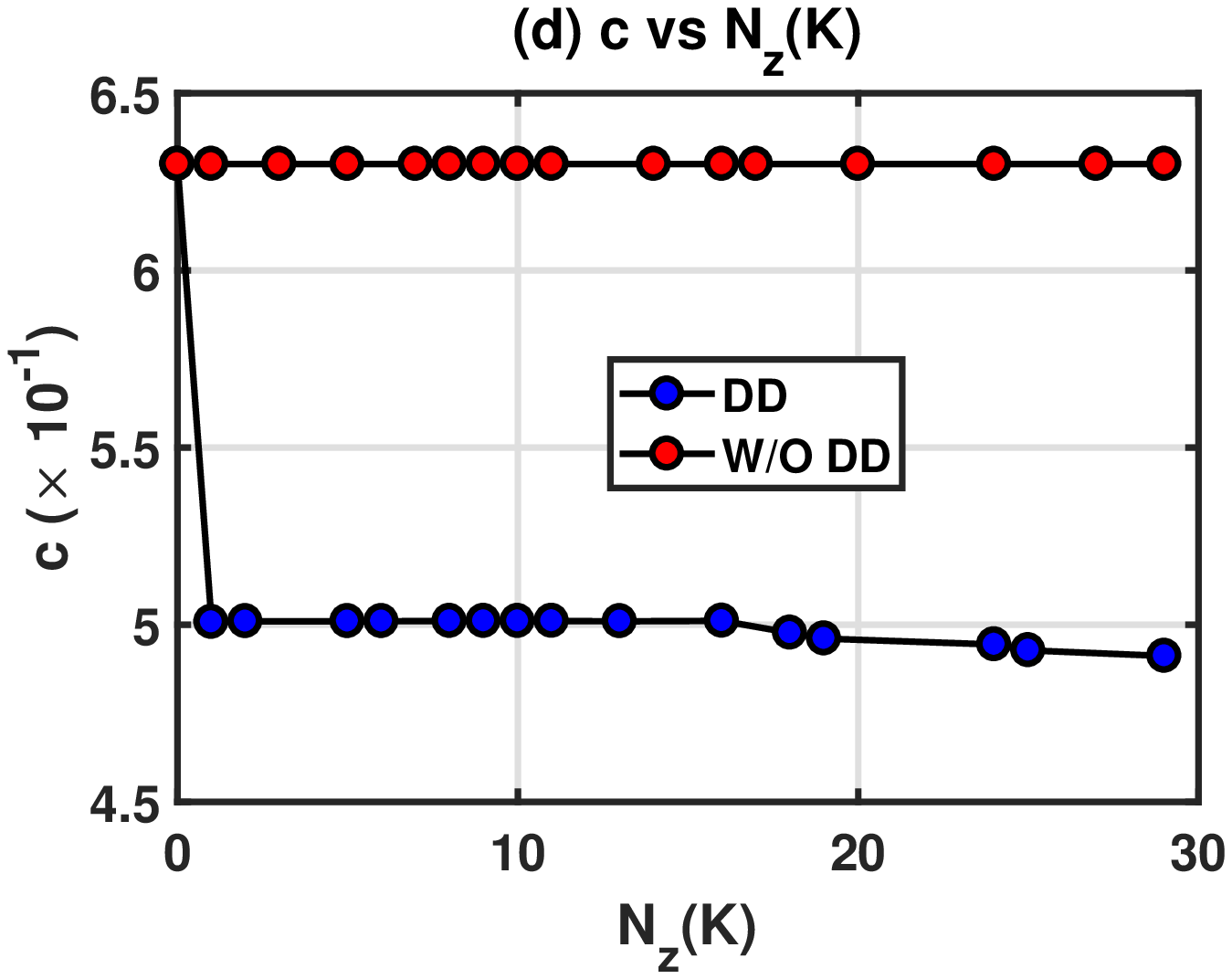}
\end{minipage}
\caption{Case B-I (a), (b), (c), (d) show $J$, $S$, $\tau_o$ and $c$ vs $N_z$ where $N_z$ is the number of zero elements of $K$. `DD' and `W/O DD' indicate Algorithm \ref{MainAlgorithm} and constant-delay algorithms respectively.}
\label{fig:Case2}
\end{figure}
\begin{figure}[hbtp]
\begin{minipage}{0.5\linewidth}
\centering
\includegraphics[scale=0.291]{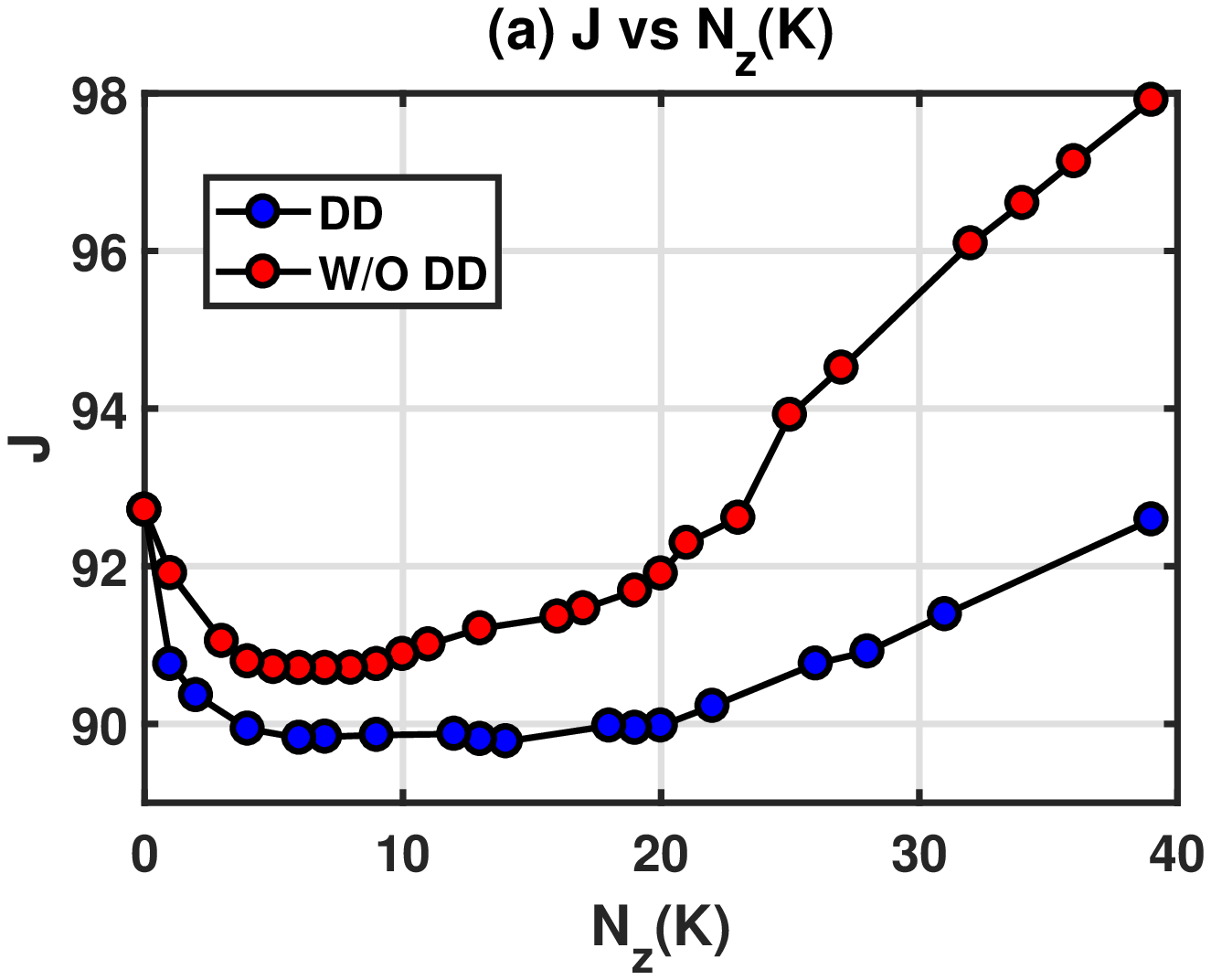} %
\end{minipage}\begin{minipage}{0.5\linewidth}
\centering
 \includegraphics[scale=0.291]{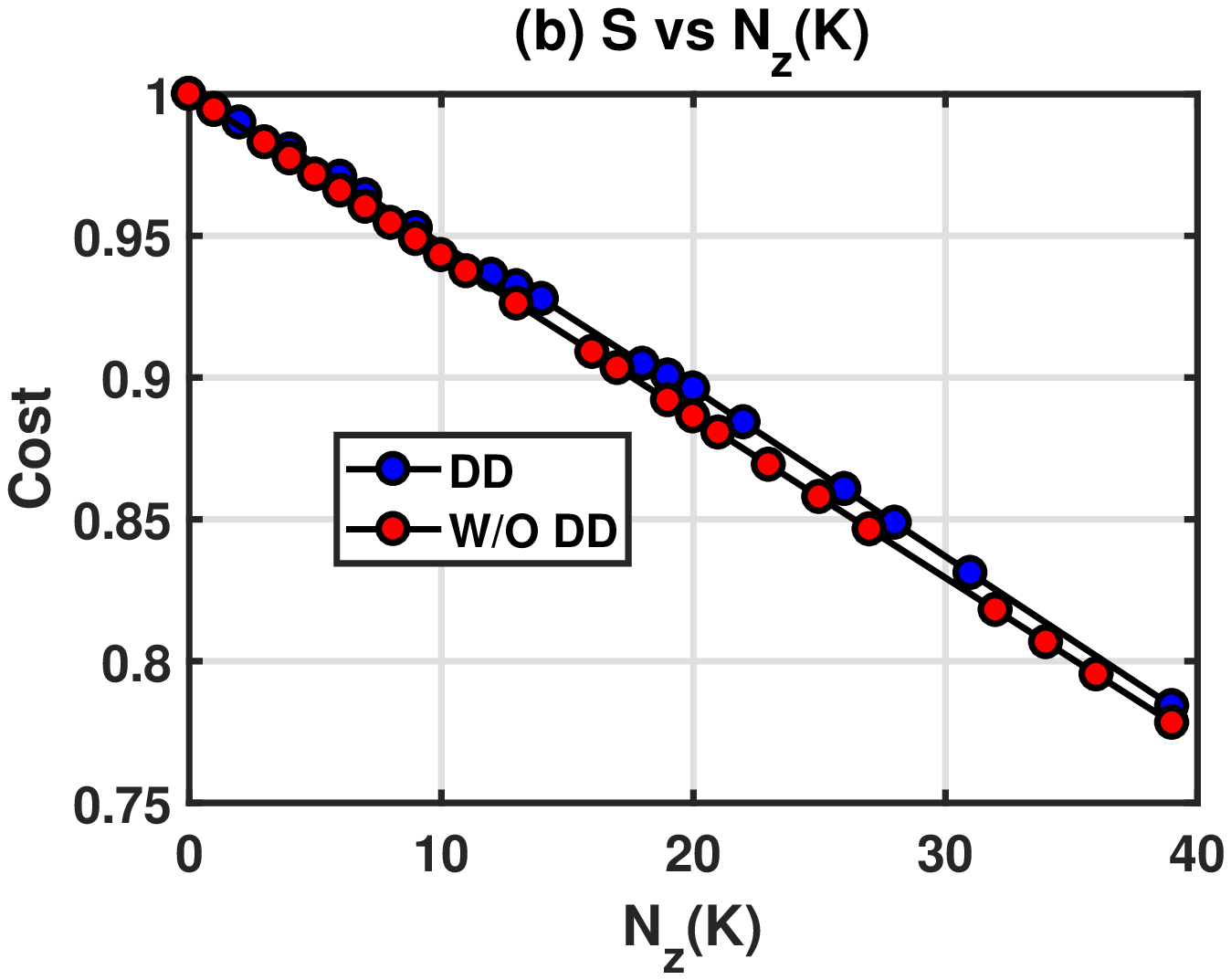}
\end{minipage}\\
\begin{minipage}{0.5\linewidth}
\centering
\includegraphics[scale=0.291]{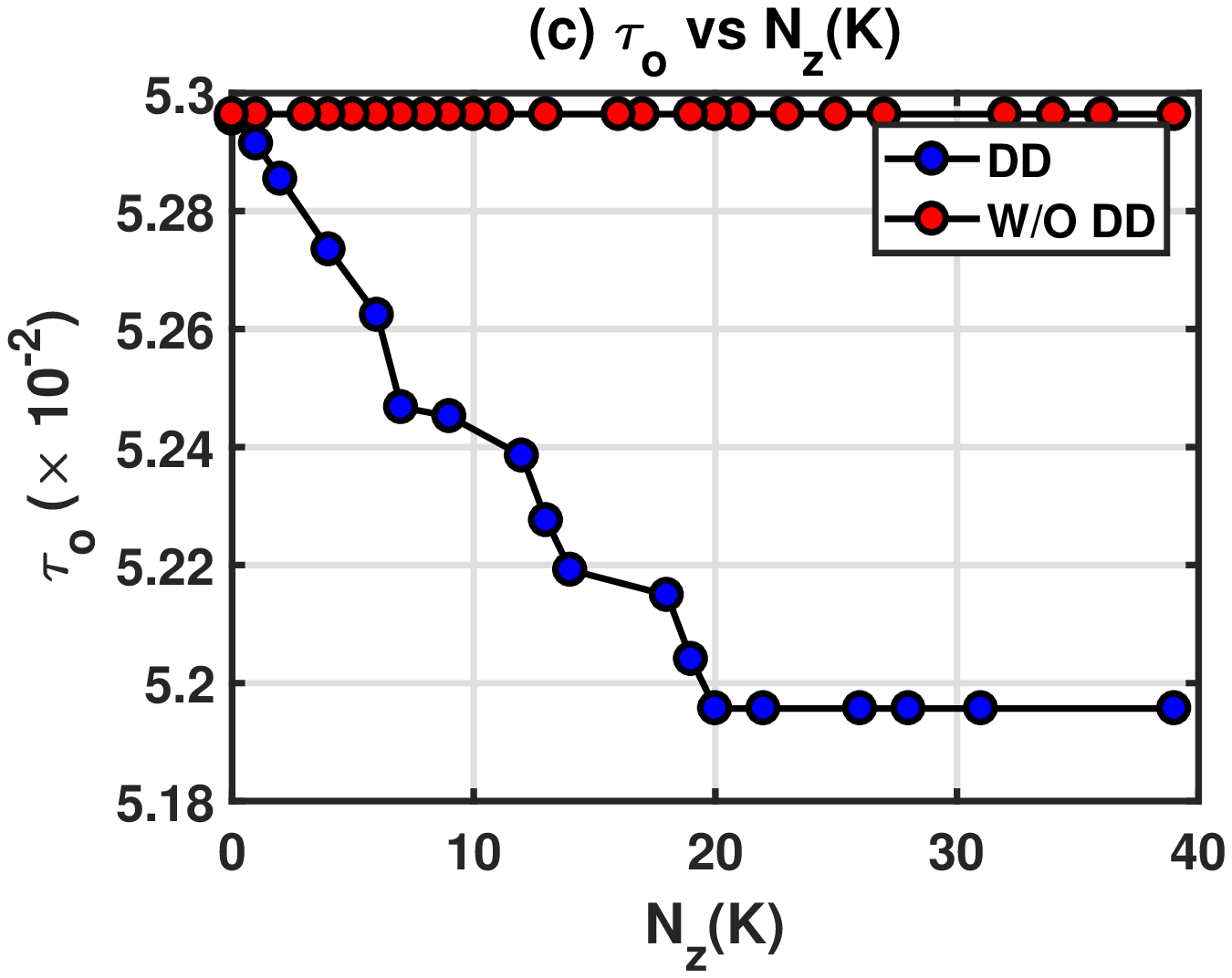} %
\end{minipage}\begin{minipage}{0.5\linewidth}
\centering
 \includegraphics[scale=0.291]{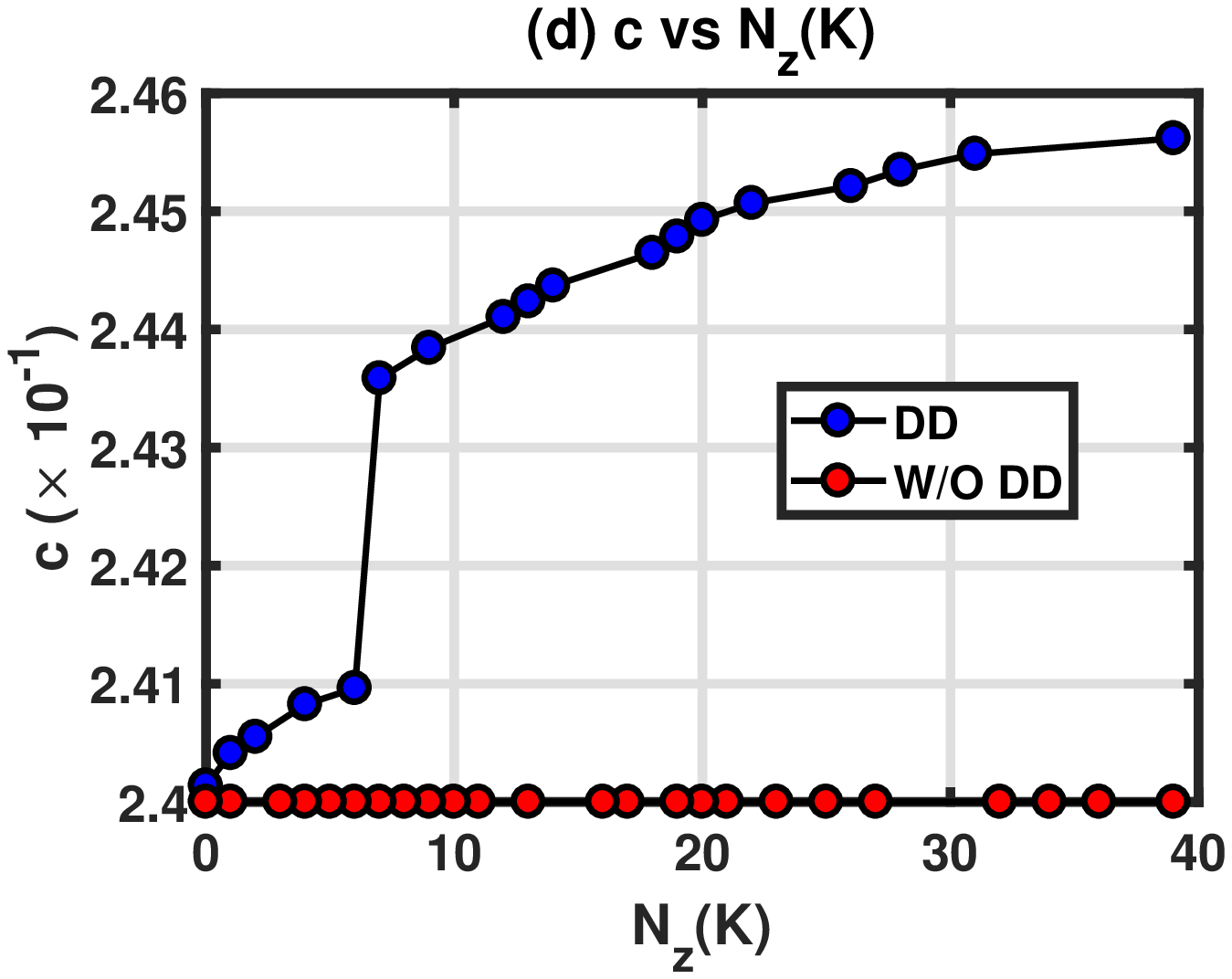}
\end{minipage}
\caption{Case B-II (a), (b), (c), (d) show $J$, $S$, $\tau_o$ and $c$ vs $N_z(K)$ where $N_z$ is the number of zero elements of $K$. `DD' and `W/O DD' indicate Algorithm \ref{MainAlgorithm} and constant-delay algorithms respectively.}
\label{fig:Case3}
\end{figure}
\subsection{Delay-Design with Bandwidth Constraints}
\label{subsec:secondsim}
\noindent We next validate Algorithm \ref{MainAlgorithm}. To illustrate its benefits, we compare it to an algorithm that consists of only the inner ADMM loop, referred to as the \textit{constant-delay} algorithm. Both algorithms start from $(K^*,\tau^*_o,\tau^*_d)\in\mc{K}$. The delays $\tau^*_o$ and $\tau^*_d$ are kept constant throughout the constant-delay algorithm. We present the simulations for two randomly generated LTI models in Case B-I and B-II with ${A}\in\mathbb{R}^{10\times 10}$ and $B=B_w=Q=R=I_n$. We denote the number of zero elements of $K$ by $N_z(K)$.
\smallskip
\par \noindent \textbf{Case B-I}: We consider $m_{cp}=53$, $m_{cc}=38$, $(c^*,\tau^*_o)=(0.63,0.053)$ for this case. Fig. \ref{fig:Case2} (a), (b), (c) and (d) show $J$, normalized bandwidth cost $S$ with respect to initial cost, $\tau_o$ and $c$, respectively. The initial delay ratio $c^*$ is such that $c_{min}=0.5007<c^{*}$, and $K^{*}$ is stabilizing for this value of $c_{min}$. As a result, Algorithm \ref{MainAlgorithm} directly moves to $c=c_{min}$ at $N_{z}(K)=1$ [see Section \ref{subsection:mainalgo}]. Further sparsification of $K$ is carried out by Algorithm \ref{MainAlgorithm} along with the simultaneous changes in $\tau_o$ and $c$, such that the bandwidth cost constraint \eqref{bw2} is satisfied, and $J$ remains optimal. Fig. \ref{fig:Case2} (c), (d) show this trade-off between $\tau_o$ and $c$, resulting in a significantly lower $S$ obtained for Algorithm \ref{MainAlgorithm} as compared to the constant-delay algorithm in Fig. \ref{fig:Case2} (b). The delay $\tau_o$ first decreases till $N_z(K)=13$, while $c$ remains nearly constant. As sparsity of $K$ increases further, $\tau_o$ begins to increase, while $c$ decreases. This indicates that the decrease in $\tau_o$ is prioritized by Algorithm \ref{MainAlgorithm} till $N_z(K)=13$. The priority later shifts to decreasing $c$ as sparsity increases. The shifting priority of one delay over another highlights the implicit relationship between $K$, $\tau_o$ and $\tau_d$.
\smallskip
\par \noindent \textbf{Case B-II}: We consider another randomly generated $A\in\mathbb{R}^{10\times 10}$ with $(c^*,\tau^*_o) =(0.24,0.089)$, $m_{cp}=21$ and $m_{cc}$ $=31$. The initial conditions result in $c_{min}=0.448>c^*$ from \eqref{constraintchange2}. However, $(K^*,c_{min})$ is an unstable tuple, and therefore, we rely on \textbf{P1}\tsb{o2} to co-design $(K,c)$. Fig. \ref{fig:Case3} (c), (d) show that as sparsity increases, Algorithm \ref{MainAlgorithm} continuously increases $c$ and decreases $\tau_o$ to maintain optimality of $J$. As shown in Proposition \ref{propositionkcktau}, a decrease in $\tau_o$ increases the bandwidth cost $S$. However, since $c$ moves towards $c_{min}$, $S$ in Fig. \ref{fig:Case3} (b) remains comparable to that of the constant-delay algorithm despite the continuous decrease in $\tau_o$. Fig. \ref{fig:Case3} (a), (b) show that as a trade-off for slightly higher $S$ from Algorithm \ref{MainAlgorithm}, we obtain a lower $J$ as compared to the constant-delay algorithm for all the sparsity levels. 
\section{Conclusion}
\par \noindent This paper presented a co-design for network delays and sparse controllers to improve the $\mc{H}_2$-performance of delayed LTI systems. The challenges of co-design arising from the implicit functional relationships between the delays, the sparse controller, and the $\mc{H}_2$-norm are overcome by developing a hierarchical algorithm, whose inner loop and outer loop are based on ADMM and SDP relaxations, respectively. Numerical simulations show the effectiveness of the design, while bringing out interesting observations about these implicit relationships. Our future work will be to extend this design to uncertain LTI models using reinforcement learning.

\section*{Appendix: Proofs}

\par \noindent \textbf{\textit{Proof of Lemma \ref{lemmatau}:}} The $ij$-th block of $\tilde{A}$ is given as:
\begin{equation}
\tilde{A}_{ij}=\begin{cases}
\sum\limits_{k=1, \ k \neq j}^N \frac{1}{\theta_j -\theta_k} I_n, \hspace{1.7cm} i=1,\ldots,N-1 \& \ i=j \\
\frac{1}{\theta_j - \theta_i}\prod\limits_{m=1, \ m\neq j, i}^N \frac{\theta_i - \theta_m}{\theta_j - \theta_m} I_n, \ \ \  i=1,\ldots,N-1 \ \& \ i\neq j,\\
A, \hspace{3.9cm} i=N, j=1,\ldots,N. 
\end{cases} \label{Aorig}
\end{equation}
Substituting \eqref{theta} above, the diagonal and off-diagonal block matrices of the first $N-1$ block rows are given by:
\begin{subequations} 
\label{lambdaAvalue}
\begin{align}
&a_{ik}=-\left({ \sin\left(\frac{(2N-i-k)\pi)}{2(N-1)} \right) \sin \left(\frac{(k-i)\pi}{2(N-1)} \right)}\right)^{-1}, \\
&\tilde{A}_{ii} =\frac{1}{\tau_o} \sum\limits_{k=1,\ k \neq i}^N  a_{ik}I_n, \ \tilde{A}_{ji} = \frac{a_{ji}}{\tau_o}\prod\limits_{m=1, \ m \neq j,i}^N \frac{a_{jm}}{a_{im}}I_n,
\end{align}
\end{subequations}
where $i,j\in \{1,\ldots,N\}$. Therefore, $\Lambda$ can be written as:
\begin{align}
 \Lambda_{ij} =\begin{cases}
 \tilde{A}_{ij}, \ \ \ \  \ i=1,\ldots,N-1, \ j=1,\ldots,N, \\
 \bb{0} , \ \ \ \ \ \ i=N, \ j=1,\ldots,N.
 \end{cases} \label{LambdaValue}
\end{align}
The proof follows from \eqref{Aorig}, \eqref{lambdaAvalue} and \eqref{LambdaValue}. $\sqbullet$
\smallskip
%%%%%%%%%%%%%%%%%%%%%%%%%%%%%%%%%%%%%%%%%%%%%%%%%
\par \noindent \textbf{\textit{Proof of Lemma \ref{lemmand}:}} { Let $\vartheta_k = \cos \left(\frac{(N-k-1)\pi}{N-1}\right)$ for $k=\{0,\ldots,N-1\}$}. From \eqref{spectral} and \eqref{acl}, $N_d$ can be written as:
\begin{align}
N_d = [l_1(-\tau_d),\ldots,l_N(-\tau_d)]^T \otimes I_n. \label{Nd}
\end{align}
Using \eqref{theta}, \eqref{lj} and $c=\nicefrac{\tau_d}{\tau_o}$ we can write
\begin{align}
l_j(-\tau_d) = \prod\limits_{m=1, \ m\neq j}^N \frac{-c-0.5(\vartheta_{m-1} -1)}{0.5(\vartheta_{j-1} - \vartheta_{m-1})} .\label{prod}
\end{align} 
Using \eqref{Nd} and \eqref{prod}, $N_d$ can be subsequently rewritten in the form of \eqref{Nddeclare} where the $j$-th row of $\Gamma$ contains the coefficients of $l_j(-\tau_d)$. From \eqref{prod}, $l_j(-\tau_d)$ is a product of $N-1$ affine terms in $c$ whose coefficients are only dependent on $N$, and hence, $\Gamma$ is constant for constant $N$. $\sqbullet$
\\
%%%%%%%%%%%%%%%%%%%%%%%%%%%%%%%%%%%%%%%%%%%%%%%
\par \noindent \textbf{\textit{Proof of Lemma \ref{lemmaunique}:}} 
The proof of uniqueness of solution of \eqref{CC} utilizes Lemma \ref{lemmatau} and \ref{lemmand}, and is similar to Theorem 2.1 in \cite{computational}. The differentiability of $P$ can be proven by utilizing the uniqueness of solution of \eqref{CC} and follows a procedure similar to Lemma 3.1 in \cite{computational}. $\sqbullet$
%%%%%%%%%%%%%%%%%%%%%%%%%%%%%%%%%%%%%%%%%
\medskip
\par \noindent \textbf{\textit{Proof of Theorem \ref{Theoremgrad}:}} The partial derivative of $J(K)$ is
\begin{equation}
 J' (K) dK = \te{Tr} ( P'(K) \mc{BB}^T) =\te{Tr}(\nabla J(K)^T dK),\label{above}
\end{equation}
where $dK\in\mathbb{R}^{m\times n}$. Post-multiplying \eqref{psiK} with $L$ and taking its trace, we obtain the following equation using \eqref{above}.
\begin{align}
\te{Tr}(dK^T\nabla J(K)) = \te{Tr} \big( (d K \circ \mc{I}_d)^T G L N_d + (d K \circ \mc{I}_o)^T G L N_o \big).\label{eq1}
\end{align}
Using the property $\te{Tr}((X\circ Y)^T Z) = \te{Tr}( X^T (Y\circ Z))$ \cite[Prob. 8.37]{schott2016matrix}, where $X,Y,Z\in\mathbb{R}^{m\times n}$ in \eqref{eq1}, we get \eqref{J(K)}. Using \eqref{taudiff}, \eqref{cdiff} in Lemma \ref{lemmaunique} and a similar procedure as above, we obtain
\begin{align}
& J'(\tau_o) = -\frac{1}{\tau_o^2} \te{Tr} ( \Lambda^T P L +L P\Lambda ), \label{33}\\
& J'(c) = \te{Tr}(N^{'}_d K_d^T G L + L G^T K_d {N^{'}_d}^T ). \label{34}
\end{align}
We can subsequently obtain \eqref{JtauJc} from \eqref{33} and \eqref{34}. $\sqbullet$
\smallskip
%%%%%%%%%%%%%%%%%%%%%%%%%%%%%%%%%%%%%%%%%%%%%%%%
\par \noindent \textbf{\textit{Proof of Theorem \ref{theoremktau}:}} Using $\phi_0$, $\phi_1$, $\psi_0$, $A_1$ and $\Delta \tilde{C}$ as stated in the theorem, we define $\phi$ and $\psi$ using Lemma \ref{lemmatau} as:
\begin{align}
& \phi = \phi_0 + \phi_1 + \phi_2, \ \phi_2=A_1^T \Delta P + \Delta P A_1^T, \nn \\
& \psi = \psi_0 + \psi_1 , \  \psi_1= \Delta\tilde{C}^T R \Delta \tilde{C}.
\end{align}
The equation $\phi+\psi =0$ is equivalent to \eqref{CC} for $(K,\omega_o,c^*)$ with $\omega_o = \nicefrac{1}{\tau_o}$ and therefore, $(K,\omega_o,c^*)$ is a stabilizing tuple if $\phi+\psi \preceq 0$ is satisfied. This inequality will be satisfied by $\phi$ and $\psi$ if they satisfy $\lambda_{max}(\phi) + \lambda_{max}(\psi) \preceq 0$ \cite[Theorem 4.3.1 (Weyl)]{matrixanalysis}. Therefore, the following inequality is a sufficient condition for stability:
\begin{align}
&\phi_0 + \phi_1+ \psi_0 + \lambda_{max} (\phi_2) + \lambda_{max}(\psi_1) \preceq 0. \label{eqbefore}
\end{align}
Equation \eqref{eqbefore} can be equivalently written as:
\begin{align}
 &\phi_0 + \phi_1+ \psi_0  +\alpha I \preceq 0, \ \alpha \geq   \lambda_{max} (\phi_2) + \lambda_{max}(\psi_1).\label{lambdamax} 
\end{align}
\par \noindent Following \cite[Theorem 4.3.50]{matrixanalysis} and \cite[Theorem 1.2]{goldberg1982numerical}, $|\lambda_{max}(\phi_2)| \leq 2 \|A^T_1 \Delta P\|$, $|\lambda_{max}(\psi_1)| = \|R^{\nicefrac{1}{2}}\Delta \tilde{C}\|^2$. Therefore, \eqref{39bc}-\eqref{39d} yield the necessary $\alpha$ for satisfying \eqref{lambdamax}. $\sqbullet$
\smallskip
%%%%%%%%%%%%%%%%%%%%%%%%%%%%%%%%%%%%%%%%%%
\par \noindent \textbf{\textit{Proof of Theorem \ref{theoremkc}:}}
Let $\phi=\sum_{i=0}^4\phi_i$ where $\phi_2=A_1\Delta L + \Delta L A^T_1$, $\phi_3=A_2 L^* + L^* A^T_2$, $\phi_4 =  A_2 \Delta L + \Delta L A^T_2$ and $A_2 = -\mc{B}\Delta K_d \Delta N^T_d$. The equation $\phi+ \mc{B}\mc{B}^T=0$ is equivalent to \eqref{LL} for $(K,\tau^*_o,c)$, and $\phi+\mc{B}\mc{B}^T\preceq 0$ implies that $(K,\tau^*_o,c)$ is a stabilizing tuple. The rest of the proof can be obtained through similar arguments as Theorem \ref{theoremktau}.
%%%%%%%%%%%%%%%%%%%%%%%%%%%%%%%%%%%%%%%%%%%%%
\medskip
\par \noindent \textbf{\textit{Derivation of \eqref{tddactual}}}: Let $\bar{k}=\texttt{vec}({K} )$. Using the property $\texttt{vec}(ABC) = (C^T \otimes A)B$, on \eqref{kmineqgrad}, we obtain the following:
\begin{align}
&(T_{dd}(\bar{k}\circ v_d) + T_{od}(\bar{v}\circ v_o) )\circ v_d + (T_{do}(\bar{k}\circ v_d) \nn\\
&+ T_{oo}(\bar{k}\circ v_o) )\circ v_o + \rho (\bar{k}\circ v_d) + \rho (\bar{k}\circ v_o) = \mu. \label{63} 
\end{align}
Since $v_d$ and $v_o$ are binary vectors, $(T_{dd}(\bar{k} \circ v_d) )\circ v_d = \big((T_{dd} \circ \hat{V}^T_d)\bar{k}\big)\circ v_d$. Furthermore, $ \big((T_{dd} \circ \hat{V}^T_d)\bar{k}\big)\circ v_d = (\hat{V}_d \circ T_{dd} \circ \hat{V}^T_d) \bar{k}$. Substituting this in \eqref{63},
\begin{align}
&(\hat{V}_d \circ T_{dd} \circ \hat{V}^T_d + \hat{V}_d \circ T_{od} \circ \hat{V}^T_o \nn\\
&+\hat{V}_o \circ T_{do} \circ \hat{V}^T_d + \hat{V}_o \circ T_{oo} \circ \hat{V}^T_o   + \rho I_{n^2}) \bar{k}  = \mu.
\end{align}
We get \eqref{tddactual} from above, thereby completing the proof. $\sqbullet$
% %%%%%%%%%%%%%%%%%%%%%%%%%%%%%%%%%%%%%%%%%%%%%%%%%%%%%%%%%%%%%

\end{document}